\documentclass[10pt]{article}
\usepackage{amsmath,amssymb,amsfonts,amsthm,epsfig,graphicx,eucal,mathrsfs}
\usepackage{latexsym,amsmath,amsfonts,graphics,subfigure}
\usepackage{accents}
\usepackage{statex}

\makeatletter
\def\widebar{\accentset{{\cc@style\underline{\mskip10mu}}}}
\def\Widebar{\accentset{{\cc@style\underline{\mskip17mu}}}}
\makeatother

\newcommand{\ignore}[1]{}

\begin{document}
\begin{center}
{\LARGE\textbf{On the skew-spectral distribution of randomly oriented graphs}}\\
\bigskip
\bigskip
Yilun Shang\\
School of Mathematical Sciences\\
Tongji University, Shanghai 200092, China\\
e-mail: \texttt{shyl@tongji.edu.cn}
\end{center}

\smallskip
\begin{abstract}

The randomly oriented graph $G_{n,p}^{\sigma}$ is an
Erd\H{o}s-R\'enyi random graph $G_{n,p}$ with a random orientation
$\sigma$, which assigns to each edge a direction so that
$G_{n,p}^{\sigma}$ becomes a directed graph. Denote by $S_n$ the
skew-adjacency matrix of $G_{n,p}^{\sigma}$. Under some mild
assumptions, it is proved in this paper that, the spectral
distribution of $S_n$ (under some normalization) converges to the
standard semicircular law almost surely as $n\rightarrow\infty$. It
is worth mentioning that our result does not require finite moments
of the entries of the underlying random matrix.

\bigskip



\end{abstract}

\bigskip
\normalsize

\section{Introduction}

Let $G$ be a simple graph with vertex set $V=\{v_1,v_2,\cdots,v_n\}$
and $G^{\sigma}$ be an oriented graph of $G$ with the orientation
$\sigma$, which assigns to each edge of $G$ a direction so that
$G^{\sigma}$ becomes a directed graph. The skew-adjacency matrix
$S(G^{\sigma})=(s_{ij})\in\mathbb{R}^{n\times n}$ is a real
skew-symmetric matrix, where $s_{ij}=1$ and $s_{ji}=-1$ if
$(v_i,v_j)$ is an arc of $G^{\sigma}$, otherwise $s_{ij}=s_{ji}=0$.
The well-known Erd\H{o}s-R\'enyi random graph model
$\mathcal{G}_{n,p}$ is a probability space \cite{1}, which consists
of all simple graphs with vertex set $V$ where each of the possible
${n\choose2}=n(n-1)/2$ edges occurs independently with probability
$p=p(n)$. For a random graph $G_{n,p}\in\mathcal{G}_{n,p}$, the
randomly oriented graph $G^{\sigma}_{n,p}$ is obtained by orienting
every edge $\{v_i,v_j\}$ $(i<j)$ in $G_{n,p}$ as $(v_i,v_j)$ with
probability $q=q(n)$ and the other way with probability $1-q$
independently of each other. Here, the superscript
$\sigma=\sigma(q)$ indicates the orientation.

The above randomly oriented graph model was first studied in
\cite{2} and a similar model based on the lattice structure (instead
of $\mathcal{G}_{n,p}$) was discussed in \cite{3}. The question of
whether the existences of directed paths between various pairs of
vertices are positively or negatively correlated has attracted some
research attention recently; see e.g. \cite{4,6,5}. Diclique
structure has been studied in \cite{20}. In this paper, we shall
explore this model from a spectral perspective. Basically, we
determine the limit spectral distribution of the random matrix
underlying the randomly oriented graph. A semicircular law
reminiscent of Wigner's famous semicircular law \cite{7} is obtained
by the moment approach (see Theorem 1 below). We mention that there
is recent increased interest in the spectral properties of oriented
graphs in classical graph theory, see e.g. \cite{10,11a,8,9}.

As is customary, we say that a graph property $\mathcal{P}$ holds
almost surely (a.s., for short) for $\mathcal{G}_{n,p}$ if the
probability that $G_{n,p}\in\mathcal{G}_{n,p}$ has the property
$\mathcal{P}$ tends to one as $n\rightarrow\infty$. We will also use
the standard Landau's asymptotic notations such as $o, O, \sim$ etc.
Let ${\bf1}_E$ be the indicator of the event $E$ and
$\textbf{i}=\sqrt{-1}$ be the imaginary unit.

\section{The results}

In this section, we characterize the spectral properties for the
skew-adjacency matrices of randomly oriented graphs.

Recall that a square matrix $M=(m_{ij})$ is said to be
skew-symmetric if $m_{ij}=-m_{ji}$ for all $i$ and $j$. It is
evident that the skew-adjacency matrix
$S_n:=S(G^{\sigma}_{n,p})=(s_{ij})\in\mathbb{R}^{n\times n}$ of the
randomly oriented graph $G^{\sigma}_{n,p}$ is a skew-symmetric
random matrix such that the upper-triangular elements $s_{ij}$
($i<j$) are i.i.d. random variables satisfying
$$
\operatorname{P}(s_{ij}=1)=pq,\quad
\operatorname{P}(s_{ij}=-1)=p(1-q)\quad\mathrm{and}\quad
\operatorname{P}(s_{ij}=0)=1-p.
$$
Hence, the eigenvalues of $S_n$ are all purely imaginary numbers. We
assume the eigenvalues are
$\textbf{i}\lambda_1,\textbf{i}\lambda_2,\cdots,\textbf{i}\lambda_n$,
where all $\lambda_i\in\mathbb{R}$.

Let $Y_n\in\mathbb{R}^{n\times n}$ be a skew-symmetric matrix whose
elements above the diagonal are 1 and those below the diagonal are
$-1$. Define a quantity
$$
r=r(p,q)=\sqrt{(1+p(1-2q))^2pq+(1-p(1-2q))^2p(1-q)}
$$
and a normalized matrix
\begin{equation}
X_n=\frac{-\textbf{i}S_n-\textbf{i}p(1-2q)Y_n}{r}.\label{1}
\end{equation}
It is straightforward to check that
$X_n=(x_{ij})\in\mathbb{C}^{n\times n}$ is a Hermitian matrix with
the diagonal elements $x_{ii}=0$ and the upper-triangular elements
$x_{ij}$ ($i<j$) being i.i.d. random variables satisfying mean
$\operatorname{E}(x_{ij})=0$ and variance
$\operatorname{Var}(x_{ij})=\operatorname{E}(x_{ij}\Widebar{x_{ij}})=1$.

In general, for a Hermitian matrix $M\in\mathbb{C}^{n\times n}$ with
eigenvalues $\mu_1(M)$, $\mu_2(M)$, $\cdots,\mu_n(M)$, the empirical
spectral distribution of $M$ is defined by
$$
F_M(x)=\frac1n\cdot\#\{\mu_i(M)|\mu_i(M)\le x,i=1,2,\cdots,n\},
$$
where $\#\{\cdots\}$ means the cardinality of a set.

\smallskip
\noindent\textbf{Theorem 1.} \itshape \quad Suppose that
$nr^2\rightarrow\infty$ and $p(1-2q)\rightarrow0$ as
$n\rightarrow\infty$. Then
$$
\lim_{n\rightarrow\infty}F_{n^{-1/2}X_n}(x)=F(x)\quad a.s.
$$
i.e., with probability 1, the empirical spectral distribution of the
matrix $n^{-1/2}X_n$ converges weakly to a distribution $F(x)$ as
$n$ tends to infinity, where $F(x)$ has the density
$$
f(x)=\frac{1}{2\pi}\sqrt{4-x^2}{\bf1}_{|x|\le2}.
$$
\normalfont

Before presenting the proof of Theorem 1, we first give a couple of
remarks.

\smallskip
\noindent\textbf{Remark 1.} The above function $F(x)$ follows the
standard semicircular distribution according to Wigner. However,
Theorem 1 extends the classical result of Wigner \cite{7}. To see
this, set $q=1/2$. The assumptions in Theorem 1 reduce to
$np\rightarrow\infty$. It is easy to check that $r=\sqrt{p}$ and
$|\operatorname{E}(x_{12}^{k+2})|=1/p^{k/2}$ if $k$ is even. Hence,
if $p=o(1)$, the condition in Wigner's semicircular law that
$\operatorname{E}(|x_{12}|^k)<\infty$ for any $k\in\mathbb{N}$ is
violated (see e.g. \cite{11,7}). In the more recent study of
spectral convergence results for Hermitian random matrices, it is
common to assume finite lower-order moments (e.g. fourth-order or
eighth-order moments) of the elements of the underlying matrices
\cite{16,18,17,15,19,14}. Therefore, our result does not fit in
these frames either.

\smallskip
\noindent\textbf{Remark 2.} Apart from Theorem 1, we can also derive
an estimate for the eigenvalues
$\textbf{i}\lambda_1,\textbf{i}\lambda_2,\cdots,\textbf{i}\lambda_n$
of the matrix $S_n$. Note that the eigenvalues of $Y_n$ are
$\mu_i(Y_n)=\textbf{i}\cot(\pi(2i-1)/2n)$ for $i=1,2,\cdots,n$. It
follows from Theorem 2.12 in \cite{12} that
$\rho(n^{-1/2}X_n)\rightarrow2$ a.s., where $\rho(\cdot)$ stands for
the spectral radius. By (\ref{1}), we have
$$
\frac{-{\bf i}S_n}{r\sqrt{n}}=\frac{X_n}{\sqrt{n}}+\frac{{\bf
i}p(1-2q)Y_n}{r\sqrt{n}}.
$$
If we arrange the eigenvalues of a Hermitian matrix
$M\in\mathbb{C}^{n\times n}$ as
$\hat{\mu}_1(M)\ge\hat{\mu}_2(M)\ge\cdots\ge\hat{\mu}_n(M)$, then
the Weyl's inequality \cite{13} implies that for all
$i=1,2,\cdots,n$,
\begin{align}
\hat{\mu}_n\left(\frac{X_n}{\sqrt{n}}\right)+\hat{\mu}_i\left(\frac{\textbf{i}p(1-2q)Y_n}{r\sqrt{n}}\right)&\le\hat{\mu}_i\left(\frac{-\textbf{i}S_n}{r\sqrt{n}}\right)\nonumber\\
&\le\hat{\mu}_1\left(\frac{X_n}{\sqrt{n}}\right)+\hat{\mu}_i\left(\frac{\textbf{i}p(1-2q)Y_n}{r\sqrt{n}}\right).\nonumber
\end{align}
Putting the above comments together, we obtain
\begin{multline}
r\sqrt{n}\left(-2+p(2q-1)\cot\left(\frac{\pi(2i-1)}{2n}\right)+o(1)\right)\le\hat{\mu}_i(-\textbf{i}S_n)\\
\le
r\sqrt{n}\left(2+p(2q-1)\cot\left(\frac{\pi(2i-1)}{2n}\right)+o(1)\right)\quad
a.s.\label{2}
\end{multline}
when $q\ge1/2$, and
\begin{multline}
r\sqrt{n}\left(-2+p(2q-1)\cot\left(\frac{\pi(2n-2i+1)}{2n}\right)+o(1)\right)\le\hat{\mu}_i(-\textbf{i}S_n)\\
\le
r\sqrt{n}\left(2+p(2q-1)\cot\left(\frac{\pi(2n-2i+1)}{2n}\right)+o(1)\right)\quad
a.s.\label{3}
\end{multline}
when $q<1/2$. Since $\hat{\mu}_i(-\textbf{i}S_n)$ is the $i$-th
largest values in the collection
$\{\lambda_1,\lambda_2,\cdots,\lambda_n\}$ by our notation, the
estimates for the eigenvalues of $S_n$ readily follow from (\ref{2})
and (\ref{3}).

Now comes the proof of Theorem 1.

\noindent\textbf{Proof of Theorem 1.} By the moment approach, it
suffices to show that the moments of the empirical spectral
distribution converge to the corresponding moments of the
semicircular law almost surely (see e.g. \cite{12}). That is,
\begin{equation}
\lim_{n\rightarrow\infty}\int x^kdF_{n^{-1/2}X_n}(x)=\int
x^kdF(x)\quad a.s.\label{4}
\end{equation}
for each $k\in\mathbb{N}$.

First note that under the assumptions of Theorem 1, it can be
checked that
\begin{equation}
\operatorname{E}(x_{12}^k)\sim\left\{\begin{array}{cc}\frac{1}{r^{k-2}}&k\equiv0(\mod
4)\\
-\frac{\textbf{i}}{r^{k-2}}&k\equiv1(\mod
4)\\
-\frac{1}{r^{k-2}}&k\equiv2(\mod
4)\\
\frac{\textbf{i}}{r^{k-2}}&k\equiv3(\mod 4)
\end{array}\right.\label{5}
\end{equation}
for any $k\in\mathbb{N}$ and $k>1$. Recall that $x_{ij}$ $(1\le
i<j\le n)$ are independently and identically distributed as
$x_{12}$, and $x_{ij}=-x_{ji}$ for all $i$ and $j$. To show
(\ref{4}), we consider the following two scenarios according to
whether $k$ is odd or even.

\textbf{(A) $k$ is odd.} Fix a $k=2t+1$ with
$t\in\mathbb{N}\cup\{0\}$. By symmetry, we have
$\int_{-2}^{2}x^kf(x)dx=0$. On the other hand, the integral on the
left-hand side of (\ref{4}) yields
\begin{eqnarray}
\int
x^kdF_{n^{-1/2}X_n}(x)&=&\frac1n\operatorname{E}\left(\operatorname{Trace}\left(\frac{X^k_n}{\sqrt{n^k}}\right)\right)=
\frac{1}{n^{1+k/2}}\operatorname{E}(\operatorname{Trace}(X_n^k))\nonumber\\
&=&\frac{1}{n^{1+k/2}}\sum_{1\le i_1,i_2,\cdots,i_k\le
n}\operatorname{E}(x_{i_1i_2}x_{i_2i_3}\cdots x_{i_ki_1}),\label{6}
\end{eqnarray}
where each summand in (\ref{6}) can be viewed as a closed walk of
length $k$ following the arcs $(v_{i_1},v_{i_2}),
(v_{i_2},v_{i_3}),\cdots,(v_{i_k},v_{i_1})$ in the complete graph
$G=K_n$ of order $n$. Clearly, each such walk contains an edge, say
$\{v_i,v_j\}$, that the total number $n_{ij}$ of times that arcs
$(v_i,v_j)$ and $(v_j,v_i)$ are traveled during this walk is odd.
Given a closed walk of length $k$, denote by $\Omega$ the set of
edges in it as described above. Thus, we have
$\Omega\not=\emptyset$. Now consider the following two cases: (A1)
there exists $\{v_i,v_j\}\in\Omega$ such that $n_{ij}=1$; (A2)
$n_{i'j'}\ge3$ for all $\{v_{i'},v_{j'}\}\in\Omega$ .

For (A1), note that the variables in the summands in (\ref{6}) are
independent and $\operatorname{E}(x_{ij})=0$. Therefore, such walks
contribute zero to the right-hand side of (\ref{6}).

For (A2), let $m$ denote the number of distinct vertices in a closed
walk of length $k$. Hence, $m$ is less than or equal to the number
of distinct vertices in a closed walk of length $2t$, in which each
edge (in either direction) appears even times. Clearly, $m\le t+1$
(the equality $m=t+1$ is attained by a walk in which each arc and
the one of opposite direction are traveled exactly once,
respectively, and all edges in the walk form a tree). Therefore,
these walks will contribute
\begin{eqnarray*}
&&\frac{1}{n^{1+k/2}}\sum_{m=1}^{t+1}\sum_{\#\{i_1,i_2,\cdots,i_k\}=m}|\operatorname{E}(x_{i_1i_2}x_{i_2i_3}\cdots
x_{i_ki_1})|\\
&\le&\frac{1}{n^{3/2+t}}\sum_{m=1}^{t+1}n^mm^k\left(\frac1r\right)^{k-2(m-1)}\\
&\le&\frac{1}{n^{3/2+t}}(t+1)n^{t+1}(t+1)^k\left(\frac1r\right)^{2t+1-2t}\\
&=&\frac{(t+1)^{k+1}}{n^{1/2}r},
\end{eqnarray*}
where the first inequality is due to (\ref{5}), (\ref{6}) and the
fact that the number of such closed walks is at most $n^mm^k$.
Consequently, combining (A1) and (A2), it follows from (\ref{6})
that
$$
\int
x^kdF_{n^{-1/2}X_n}(x)=O\left(\frac{1}{n^{1/2}r}\right)\rightarrow0
$$
as $n\rightarrow\infty$, by our assumptions. We complete the proof
of (\ref{4}) for odd $k$.

\textbf{(B) $k$ is even.} Fix a $k=2t$ with
$t\in\mathbb{N}\cup\{0\}$. We have
\begin{eqnarray}
\int_{-2}^{2}x^kf(x)dx&=&\frac{1}{2\pi}\int_{-2}^2x^k\sqrt{4-x^2}dx=\frac{1}{\pi}\int_0^2x^{2t}\sqrt{4-x^2}dx\nonumber\\
&=&\frac{2^{2t+1}}{\pi}\int_0^1y^{t-1/2}(1-y)^{1/2}dy\nonumber\\
&=&\frac{2^{2t+1}}{\pi}\cdot\frac{\Gamma(t+1/2)\Gamma(3/2)}{\Gamma(t+2)}\nonumber\\
&=&\frac{1}{t+1}{2t\choose t}.\label{7}
\end{eqnarray}

Given a closed walk of length $k$ in $K_n$, we still set $m$ as the
number of distinct vertices in it. To analyze the terms in
(\ref{6}), we consider the following three cases: (B1) there exists
an edge, say $\{v_i,v_j\}$, in the closed walk such that the total
number of times that arcs $(v_i,v_j)$ and $(v_j,v_i)$ are traveled
during this walk is odd; (B2) no such $\{v_i,v_j\}$ exists, and
$m\le t$; (B3) no such $\{v_i,v_j\}$ exists, and $m=t+1$. Note that
if each edge (in either direction) of the closed walk appears even
times, we have $m\le t+1$. The equality holds if and only if each
arc and the one of opposite direction are traveled exactly once,
respectively, and all edges in the walk form a tree.

For (B1), we argue similarly as in (A1) and know that the
contribution to the right-hand side of (\ref{6}) is zero.

For (B2), an analogous derivation as in (A2) reveals that the
contribution to the right-hand side of (\ref{6}) amounts to
\begin{eqnarray*}
&&\frac{1}{n^{1+k/2}}\sum_{m=1}^t\sum_{\#\{i_1,i_2,\cdots,i_k\}=m}|\operatorname{E}(x_{i_1i_2}x_{i_2i_3}\cdots
x_{i_ki_1})|\\
&\le&\frac{1}{n^{1+t}}\sum_{m=1}^tn^mm^k\left(\frac1r\right)^{k-2(m-1)}\\
&\le&\frac{1}{n^{1+t}}\cdot t\cdot n^t\cdot t^k\cdot\left(\frac1r\right)^{2t-2(t-1)}\\
&=&\frac{t^{k+1}}{nr^2}.
\end{eqnarray*}

For (B3), noting that
$\operatorname{E}(x_{12}x_{21})=-\operatorname{E}(x_{12}^2)=1$ and
the independence of the variables, we obtain that each term
$\operatorname{E}(x_{i_1i_2}x_{i_2i_3}\cdots x_{i_ki_1})$ in
(\ref{6}) equals 1. Recall that a combinatorial result \cite[Lemma
2.4]{18} says that the number of the closed walks of length $2t$ on
$t+1$ vertices, which satisfy that each each arc and the one of
opposite direction both appear exactly once, and all edges in the
walk form a tree, equals $\frac{1}{t+1}{2t\choose t}(t+1)!$. Since
there are ${n\choose t+1}$ choices of a set of $t+1$ vertices, we
conclude that the contribution to the left-hand side of (\ref{6})
amounts to
$$
\frac{1}{n^{1+k/2}}\cdot\frac{1}{t+1}{2t\choose
t}(t+1)!\cdot{n\choose
t+1}=\frac{n(n-1)\cdots(n-t)}{n^{1+t}}\cdot\frac{1}{t+1}{2t\choose
t}.
$$

Finally, combining (B1), (B2) and (B3), it follows from (\ref{6})
that
\begin{eqnarray*}
\int
x^kdF_{n^{-1/2}X_n}(x)&=&O\left(\frac{1}{nr^2}\right)+\frac{n(n-1)\cdots(n-t)}{n^{1+t}}\cdot\frac{1}{t+1}{2t\choose
t}\\
&\rightarrow&\frac{1}{t+1}{2t\choose t},
\end{eqnarray*}
as $n\rightarrow\infty$, by our assumptions. In view of (\ref{7}),
the proof of (\ref{4}) for even $k$ is complete. $\Box$

To conclude the paper, we simulate the randomly oriented graph model
and computed the eigenvalue distribution for the matrix
$n^{-1/2}X_n$ (see Figure 1). The simulation results show a perfect
agreement with our theoretical prediction. For future work, it would
be interesting to explore some other properties (see e.g.
\cite{bb,21}) in the setting of randomly oriented graphs.

\begin{figure}[htb] \centering
{\subfigure[]{\includegraphics[width=.49\textwidth]{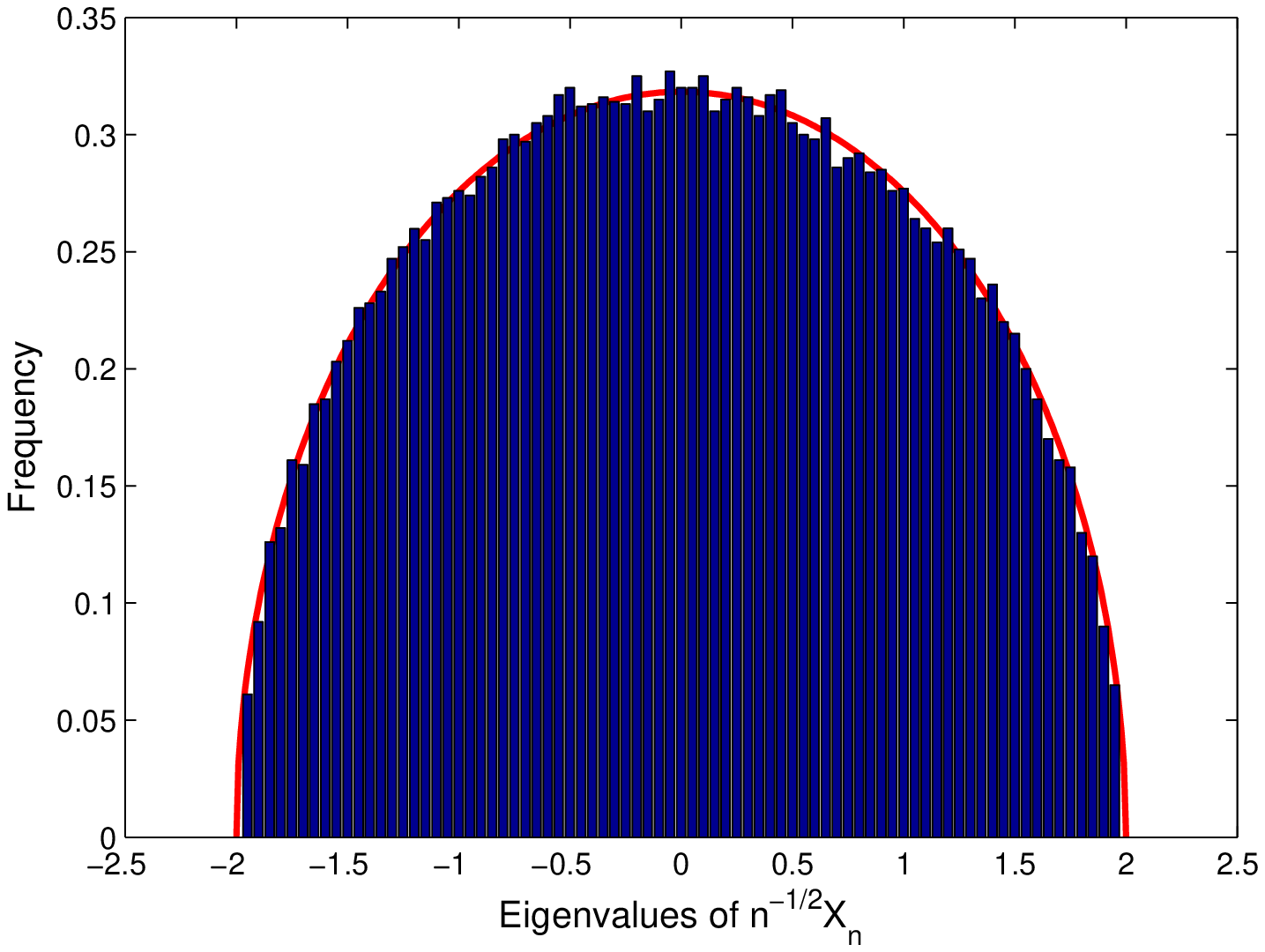}}}
{\subfigure[]{\includegraphics[width=.49\textwidth]{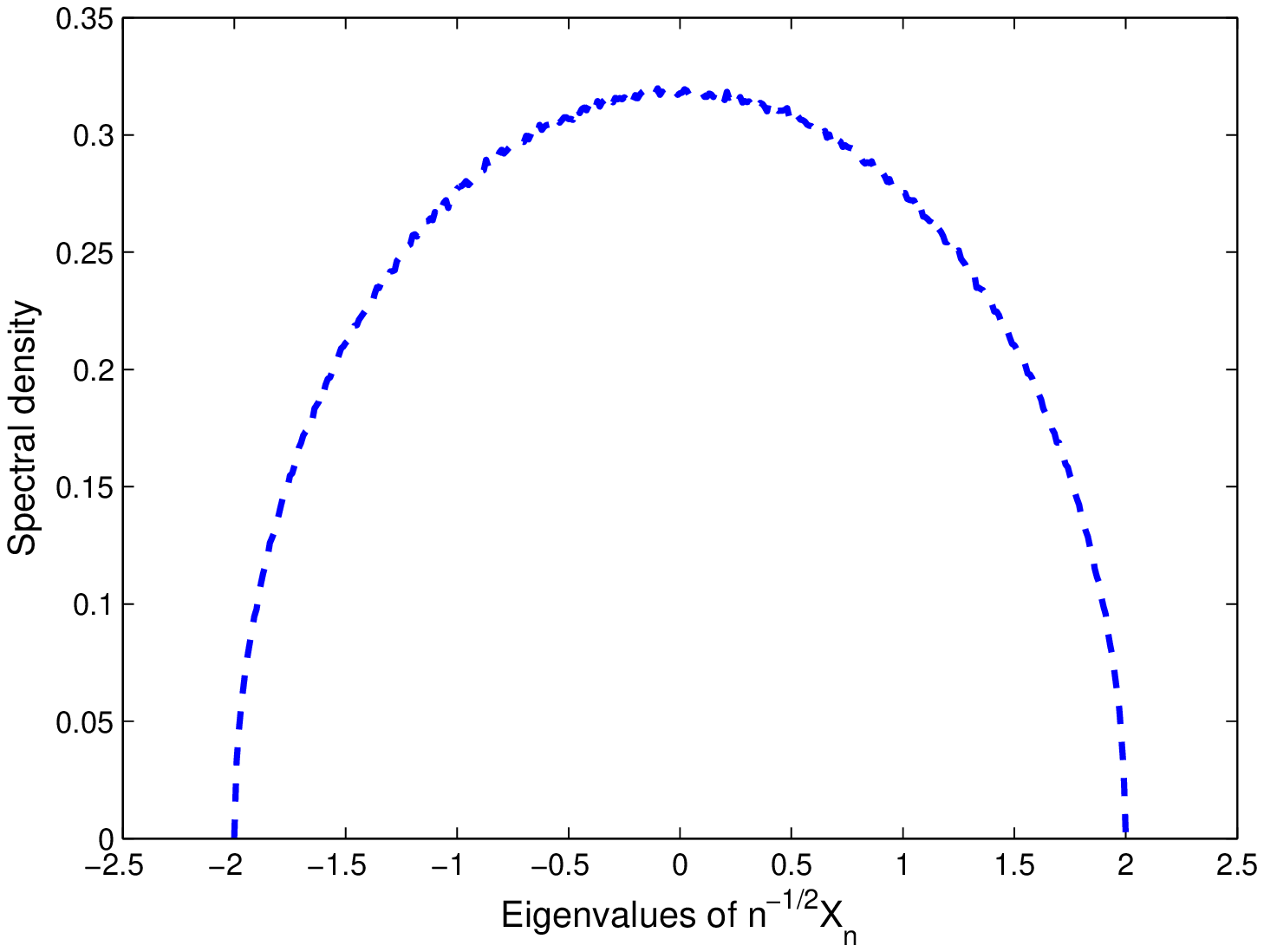}}}
\caption{Limiting skew-spectral distribution for
$G_{n,p}^{\sigma(q)}$ with $n=1000$, $p=0.1$ and $q=0.5$. (a)
Histogram of the spectrum of $n^{-1/2}X_n$. A solid line shows the
semicircular distribution for comparison. (b) The average spectral
density for $n^{-1/2}X_n$ over 500 graphs.}
\end{figure}

\section*{Acknowledgements}

The author is indebted to the referees for helpful comments. The
work is supported by the National Natural Science Foundation of
China (11505127) and the Shanghai Pujiang Program (15PJ1408300).

\end{document}